\documentclass[a4,12pt]{article}
\usepackage{amssymb}
\usepackage{amstext}
\title{The joint distribution of occupation times of skip-free Markov processes and 
a class of multivariate exponential distributions}
\author{Kshitij Khare}

\begin{document}
\maketitle

\begin{abstract}
For a skip-free Markov process on $\mathbb{Z}^+$ with generator matrix $Q$, we 
evaluate the joint Laplace transform of the occupation times before hitting the 
state $n$ (starting at $0$). This Laplace transform has a very straightforward and 
familiar expression. We investigate the properties of this Laplace transform, 
especially the conditions under which the occupation times form a Markov chain. 
\end{abstract}

\section{Introduction}
Consider a process $\{X_t\}_{t \geq 0}$ on $\{0,1,2, ..........\}$ with generator 
$Q := ((q_{ij}))_{0 \leq i,j < \infty}$ satisfying 
\begin{equation} \label{asmpnskpfr}
q_{i,i+1} > 0, q_{ij} = 0 \mbox{ if } j > i + 1 \mbox{ and } \sum_{j=0}^{i+1} 
q_{ij} = 0 \mbox{ for all } i \geq 0. 
\end{equation}

\noindent
Such processes are called skip-free (to the left) processes. Note that birth and 
death processes are special cases of the process described above. Let $T_n$ be the 
first time when the process reaches $n$. As usual, $l_{T_n}^i$ denotes the amount 
of time $\{X_t\}_{t \geq 0}$ spends in state $i$ before it reaches $n$. Let $Q_n$ 
denote the upper $n \times n$ submatrix of $Q$. We are interested in the joint 
distribution of $\left\{ l_{T_n}^i \right\}_{0 \leq i \leq n-1}$ starting at $0$. 
We prove that 
\begin{equation} \label{jtlpctrnfm}
{\bf{E}}_0 \left[ e^{-\sum_{i=0}^{n-1} d_i l_{T_n}^i} \right] = \frac{|-Q_n|}
{|D - Q_n|}, 
\end{equation}

\noindent
where $d_i \geq 0, \; i = 0,1, ... ,n-1$ and $D$ is a diagonal matrix with diagonal 
entries $d_i, \; i = 0,1, ... ,n-1$. This formula generalizes the results of Kent 
\cite{kentmedocc} about birth and death processes. Using the identity 
(\ref{jtlpctrnfm}) we prove that $\left\{ l_{T_n}^i \right\}_{0 \leq i \leq n-1}$ is 
a Markov chain iff $Q_n$ is a tridiagonal matrix. For these purposes, we will 
derive an identity for the joint Laplace transfrom of occupation times for a 
general class of Markov processes and then specialize to skip-free processes. 

\section{An identity for a general class of Markov processes}
Let $\{X_t\}_{t \geq 0}$ be a Markov process with finite state space $\mathcal{X}$. 
Let $Q = ((q_{xy}))_{x,y \in \mathcal{X}}$ be the corresponding generator matrix. 

\indent
Let $S_i$ be the random time corresponding to the $i^{th}$ transition for the 
Markov process $\{X_t\}_{t \geq 0}, \; i = 1,2,3, ..........$ Set $S_0 = 0$. Let 
$\{Y_i\}_{i \geq 0} := \{X_{S_i}\}_{i \geq 0}$ be the embedded discrete-time 
Markov chain with one step transition probabilities 
$$
p_{xy} = \frac{-q_{xy}}{q_{xx}} \mbox{ for } x \neq y, \; p_{xx} = 0. 
$$

\noindent
Let $P := ((p_{xy}))_{x,y \in \mathcal{X}}$ and let $Q^{diag}$ denote the 
diagonal matrix with diagonal entries same as $Q$. Then, 
\begin{equation} \label{idntyembed}
Q = Q^{diag} (I - P). 
\end{equation}

\noindent
As usual, let 
\begin{eqnarray*}
l_t^x := \int_0^t 1_{\{X_s = x\}} ds, \; t \geq 0, x \in \mathcal{X} 
\end{eqnarray*}

\noindent
denote the occupation time of the Markov process $\{X_t\}_{t \geq 0}$ in the state 
$x$ till time $t$. 

\indent
We describe a typical path of $\{X_t\}_{t \geq 0}$. The process starts at an 
initial state $Y_0$. The process stays at $Y_n$ during $[S_n, S_{n+1})$ for 
$n = 0,1,2, ..........$ and at a random time $\xi = S_{\eta + 1}$ jumps from the 
state $Y_{\eta}$ to a ``cemetery" $\Delta$ (not included in $\mathcal{X}$), and 
stays there forever. (Here $\xi$ takes non-negative real values and $\eta$ takes 
non-neagtive integer values). The analysis presented below requires the 
assumptions that 

\smallskip

\noindent
(i) $Q^{-1}$ exists. 

\noindent
(ii) ${\bf{P}}_x (\eta < \infty) = 1 \; \forall x \in \mathcal{X}$. 

\smallskip

\noindent
Let $\mathbf{d} = \{d_u\}_{u \in \mathcal{X}}$  be arbitrary with $d_u \geq 0$. 
Let $D$ denote the diagonal matrix with diagonal entries $\{d_u\}_{u \in 
\mathcal{X}}$. Note that $-(D - Q)$ is the generator of a Markov process 
$\{\bar{X}_t\}_{t \geq 0}$ with the same structure as $\{X_t\}_{t \geq 0}$ except 
that at every state $x \in \mathcal{X}$, there is an additional killing rate of 
$d_x$. Let $\{\bar{Y}_m\}_{m \geq 0}$ be the embedded discrete-time Markov chain 
corresponding to $\{\bar{X}_t\}_{t \geq 0}$. Let us establish the change of 
measure formula from $\{\bar{Y}_m\}_{m \geq 0}$ to $\{Y_m\}_{m \geq 0}$. 
\begin{eqnarray*}
& & {\bf{P}}_x \{\bar{Y}_1 = y_1, \bar{Y}_2 = y_2, ... , \bar{Y}_n = y_n\}\\
&=& \prod_{i=1}^n \frac{q_{y_{i-1}y_i}}{-q_{y_{i-1}y_{i-1}} + d_{y_{i-1}}} \; \; 
(\mbox{with } y_0 = x)\\
&=& \left( \prod_{i=1}^n \frac{-q_{y_{i-1}y_{i-1}}}{-q_{y_{i-1}y_{i-1}} + 
d_{y_{i-1}}} \right) {\bf{P}}_x \{Y_1 = y_1, Y_2 = y_2, ... , Y_n = y_n\}
\end{eqnarray*}

\noindent
Hence, 
\begin{equation} \label{chgmridnty}
{\bf{E}}_x [F(\bar{Y}_1, \bar{Y}_2, ... , \bar{Y}_n)] = {\bf{E}}_x \left[ 
\left( \prod_{i=1}^n \frac{-q_{Y_{i-1}Y_{i-1}}}{-q_{Y_{i-1}Y_{i-1}} + 
d_{Y_{i-1}}} \right) F(Y_1, Y_2, ... , Y_n) \right] 
\end{equation}

\noindent
for each bounded Borel measurable function $F$. 

\smallskip

\noindent
We evaluate ${\bf{E}}_x \left[ e^{-\sum_{u \in \mathcal{X}} d_u l_{\infty}^u} 
\right]$ which is the joint Laplace transform of the occupation times 
$\{l_{\infty}^u\}_{u \in \mathcal{X}}$ evaluated at $\mathbf{d} = \{d_u\}_{u \in 
\mathcal{X}}$. The argument is inspired by Dynkin \cite{dynkmrkvps}. 

\noindent
Note that, 
\begin{eqnarray*}
\sum_{u \in \mathcal{X}} d_u l_{\infty}^u = \sum_{i=0}^{\eta} d_{Y_i} 
(S_{i+1} - S_i). 
\end{eqnarray*}

\noindent
Hence, 
\begin{eqnarray*}
& & {\bf{E}}_x \left[ e^{-\sum_{u \in \mathcal{X}} d_u l_{\infty}^u} \right]\\
&=& \sum_{n=0}^{\infty} {\bf{E}}_x \left[ e^{-\sum_{i=0}^n d_{Y_i} 
(S_{i+1} - S_i)} 1_{\{\eta = n\}} \right]\\
&=& \sum_{n=0}^{\infty} {\bf{E}}_x \left[ e^{-\sum_{i=0}^n d_{Y_i} 
(S_{i+1} - S_i)} \sum_{y \in \mathcal{X}} 1_{\{Y_n = y, Y_{n+1} = \Delta\}} 
\right]\\
&=& \sum_{y \in \mathcal{X}} \sum_{n=0}^{\infty} {\bf{E}}_x \left[ 
e^{-\sum_{i=0}^n d_{Y_i} (S_{i+1} - S_i)} 1_{\{Y_n = y, Y_{n+1} = \Delta\}} 
\right]\\
&=& \sum_{y \in \mathcal{X}} \sum_{n=0}^{\infty} {\bf{E}}_x \left[ {\bf{E}}_x 
\left[ e^{-\sum_{i=0}^n d_{Y_i} (S_{i+1} - S_i)} \mid \{Y_m\}_{m \geq 0} 
\right] 1_{\{Y_n = y, Y_{n+1} = \Delta\}} \right]\\
&=& \sum_{y \in \mathcal{X}} \sum_{n=0}^{\infty} {\bf{E}}_x 
\left[ \left( \prod_{i=0}^n \frac{-q_{Y_i Y_i}}{-q_{Y_i Y_i} + d_{Y_i}} \right) 
1_{\{Y_n = y, Y_{n+1} = \Delta\}} \right]
\end{eqnarray*}

\noindent
The previous equality follows from the fact that conditioned on 
$\{Y_m\}_{m \geq 0}$, the intermediate transition times 
$\{S_{i+1} - S_i\}_{i \geq 0}$ are independent and have 
$Exponential(-q_{Y_i Y_i})$ distribution for $i = 0,1,2, ..........$ 

\noindent
By Markov property and (\ref{chgmridnty}), 
\begin{eqnarray*}
& & {\bf{E}}_x \left[ e^{-\sum_{u \in \mathcal{X}} d_u l_{\infty}^u} \right]\\
&=& \sum_{y \in \mathcal{X}} \sum_{n=0}^{\infty} \frac{-q_{yy}}{-q_{yy} + d_y} 
{\bf{E}}_x \left[ \left( \prod_{i=1}^n \frac{-q_{Y_{i-1}Y_{i-1}}}
{-q_{Y_{i-1}Y_{i-1}} + d_{Y_{i-1}}} \right) 1_{\{Y_n = y\}} \right] 
{\bf{P}}(Y_1 = \Delta \mid Y_0 =y)\\
&=& \sum_{y \in \mathcal{X}} \sum_{n=0}^{\infty} \frac{-q_{yy}}{-q_{yy} + d_y} 
{\bf{E}}_x \left[ 1_{\{\bar{Y}_n = y\}} \right] p(y, \Delta)\\
&=& \sum_{y \in \mathcal{X}} -q_{yy} {\bf{E}}_x \left[ \int_0^{\infty} 
1_{\{\bar{X}_s = y\}} ds \right] p(y, \Delta)\\
&=& \sum_{y \in \mathcal{X}} -q_{yy} (D - Q)^{-1} (x,y) \; p(y, \Delta)\\
&=& \frac{|-Q|}{|D - Q|} \left( \sum_{y \in \mathcal{X}} -q_{yy} 
\frac{|(D - Q)^{(x,y)}|}{|-Q|} p(y, \Delta ) \right)
\end{eqnarray*}

\noindent
Here, $|(D - Q)^{(x,y)}|$ represents the determinant of the matrix obtained after 
removing the $y^{th}$ row and the $x^{th}$ column from $D - Q$, multiplied by $-1$ 
if the $x^{th}$ row and the $y^{th}$ row of $D - Q$ differ by an odd number of 
rows. Thus, the joint Laplace transform of $\{l_{\infty}^u\}_{u \in \mathcal{X}}$ 
is given by 
\begin{equation} \label{lpctrocctm}
{\bf{E}}_x \left[ e^{-\sum_{u \in \mathcal{X}} d_u l_{\infty}^u} \right] = 
\frac{|-Q|}{|D - Q|} \left( \sum_{y \in \mathcal{X}} \frac{|(D - Q)^{(x,y)}|}
{|-Q|} \left( -\sum_{u \in \mathcal{X}} q_{yu} \right) \right). 
\end{equation}

\noindent
In \cite{eikamarosh}, (\ref{lpctrocctm}) is derived using Kac's moment formula. 
Note that we did not require any reversibility assumption on $Q$ in the above 
argument. 

\section{Application to skip-free processes}
Let us investigate the case of skip-free processes with generator $Q$ given by 
(\ref{asmpnskpfr}). We are interested in the joint distribution of 
$\left\{ l_{T_n}^i \right\}_{0 \leq i \leq n-1}$ starting at state $0$. It follows 
by the structure of $Q$ that the joint distribution of $\left\{ l_{T_n}^i 
\right\}_{0 \leq i \leq n-1}$ starting at $0$ is the same as the joint 
distribution of the infinite occupation times of a Markov chain with generator 
matrix $Q_n$, starting at $0$. Hence, by (\ref{lpctrocctm}), for arbitrary 
$d_i \geq 0, \; i = 0,1, ... , n-1$, 
\begin{eqnarray*}
{\bf{E}}_0 \left[ e^{-\sum_{i=0}^{n-1} d_i l_{T_n}^i} \right] 
&=& \frac{|-Q_n|}{|D - Q_n|} \left( \sum_{j=0}^{n-1} \frac{|(D - Q_n)^{(0,j)}|}
{|-Q_n|} \left( -\sum_{i=0}^{n-1} q_{ji} \right) \right)\\
&=& -\frac{|-Q_n|}{|D - Q_n|} \frac{|(D - Q_n)^{(0,n-1)}|}{|-Q_n|} 
\sum_{i=0}^{n-1} q_{n-1,i} 
\end{eqnarray*}

\noindent
The previous equality follows from the fact that 
\begin{eqnarray*}
\sum_{i=0}^{n-1} q_{ji} = 0 \mbox{ for } j < n - 1. 
\end{eqnarray*}

\noindent
Note that if we remove the $(n-1)^{th}$ row and $0^{th}$ column of $D - Q_n$, the 
resulting matrix is upper triangular with diagonal entries not depending on 
$d_i, \; i = 0,1, ... , n-1$. Hence, 
\begin{eqnarray*}
{\bf{E}}_0 \left[ e^{-\sum_{i=0}^{n-1} d_i l_{T_n}^i} \right] \propto 
\frac{1}{|D - Q_n|}. 
\end{eqnarray*}

\noindent
By substituting  $d_i = 0$ for $i = 0,1, ... , n-1$, we get, 
\begin{eqnarray*}
{\bf{E}}_0 \left[ e^{-\sum_{i=0}^{n-1} d_i l_{T_n}^i} \right] = 
\frac{|-Q_n|}{|D - Q_n|}. 
\end{eqnarray*}

\noindent
Let us look more closely at the joint Laplace transform in (\ref{jtlpctrnfm}). Note 
that, if $Q_n^{-1} := ((q_n^{ij}))_{0 \leq i,j \leq n-1}$, then 
$$
{\bf{E}}_0 \left[ e^{-d_i l_{T_n}^i} \right] = \frac{-\frac{1}{q_n^{ii}}}{d_i - 
\frac{1}{q_n^{ii}}} \mbox{ for } d_i \geq 0, \; 0 \leq i \leq n-1. 
$$

\noindent
Hence, the marginal distribution of $l_{T_n}^i$ is $Exponential \left( \frac{1}
{q_n^{ii}} \right)$ for $0 \leq i \leq n-1$, which means the joint distribution of 
$\left\{ l_{T_n}^i \right\}_{0 \leq i \leq n-1}$ is a multivariate exponential 
distribution. 

\indent
It is well known that if $\mathbf{\eta}, \tilde{\mathbf{\eta}} \stackrel{i.i.d}
{\sim} MVN_n (\mathbf{0}, \Sigma)$, i.e. $\mathbf{\eta} = (\eta_1, \eta_2, ... ,
\eta_n)$ and $\tilde{\mathbf{\eta}} = (\tilde{\eta}_1, \tilde{\eta}_2, ... , 
\tilde{\eta}_n)$ are multivariate normal with mean $\mathbf{0}$ and covariance 
matrix $\Sigma$, then 
\begin{eqnarray*}
{\bf{E}} \left[ e^{-\sum_{i=1}^n \frac{d_i}{2} (\eta_i^2 + \tilde{\eta}_i^2)} 
\right] = \frac{\left| \Sigma^{-1} \right|}{\left| D + \Sigma^{-1} \right|} 
\mbox{ for } d_i \geq 0, \; i = 1,2, ... ,n. 
\end{eqnarray*}

\noindent
Note that if $M = \tilde{D} \Sigma^{-1} \tilde{D}^{-1}$, where $\tilde{D}$ is a 
diagonal matrix, then 
$$
\frac{\left| \Sigma^{-1} \right|}{\left| D + \Sigma^{-1} \right|} = \frac{|M|}
{|D + M|}. 
$$

\noindent
If $Q$ is reversible with respect to a probability measure $\pi$, it is possible to 
symmetrize $Q$ by diagonal conjugation, and hence we can identify the joint 
distribution of $\left\{ l_{T_n}^i \right\}_{0 \leq i \leq n-1}$ with the sum of 
componentwise squares of two independent $MVN_{n} \left( \mathbf{0}, (-Q_n^*)^{-1} 
\right)$ random vectors multiplied by $\frac{1}{2}$. $Q_n^*$ is given by 
$$
Q_n^* (i,j) = q_{ij} \sqrt{\frac{\pi(i)}{\pi(j)}}, \; 0 \leq i,j \leq n-1. 
$$

\noindent
But a generator for a skip-free process is diagonally conjugate to a symmetric 
matrix if and only if the process is a birth and death process. We refer to 
skip-free processes which are not birth and death processes as {\it{strictly}} 
skip-free processes. Hence, if $Q_n$ is not tridiagonal, then the joint 
distribution of $\left\{ l_{T_n}^i \right\}_{0 \leq i \leq n-1}$ can not be 
identified with componentwise sums of squares of independent multivariate normal 
random vectors multiplied by $\frac{1}{2}$. 

\section{A complex Gaussian measure}
If $A$ is real non-symmetric and positive definite, we still can interpret 
$$
\phi_A (\mathbf{d}) = \frac{|A|}{|D + A|} \mbox{ for } \mathbf{d} = 
(d_1, d_2, ... ,d_n) \mbox{ with } d_i \geq 0, \; i = 1,2, ... ,n, 
$$

\noindent
as the Laplace transform of a signed measure. This will be proved using the methods 
in \cite{lejandknim}. 

\indent
Let $C := \frac{A + A^T}{2}$ and $B := \frac{A - A^T}{2}$. Clearly, $C$ is 
symmetric, $B$ is skew-symmetric and $A = C + B$. Define the measure $\mu_A$ on 
$\mathbb{C}^n$ by 
$$
d \mu_A (\mathbf{z}) = \left( e^{-\frac{\mathbf{z}^T A \bar{\mathbf{z}}}{2}} + 
e^{-\frac{\bar{\mathbf{z}}^T A \mathbf{z}}{2}} \right) d \mathbf{z}, 
$$

\noindent
where $\mathbf{z} = \mathbf{x} + i \mathbf{y}$ and $d \mathbf{z} = \prod_{j=1}^n 
dx_j \prod_{j=1}^n dy_j$. Note that $\mathbf{z}^T A \bar{\mathbf{z}}$ and 
$\bar{\mathbf{z}}^T A \mathbf{z}$ are complex conjugates, hence $\mu_A$ is a 
real valued (although possibly signed) measure. Also, the transformation 
$\mathbf{y} \rightarrow -\mathbf{y}$ (where $\mathbf{z} = \mathbf{x} + i 
\mathbf{y}$) gives, 
$$
\int_{\mathbb{C}^n} e^{-\frac{\mathbf{z}^T A \bar{\mathbf{z}}}{2}} d \mathbf{z} = 
\int_{\mathbb{C}^n} e^{-\frac{\bar{\mathbf{z}}^T A \mathbf{z}}{2}} d \mathbf{z}. 
$$

\noindent
To avoid confusion, we clarify that $\int_{\mathbb{C}^n} f(\mathbf{z}) d 
\mathbf{z}$ stands for $\int_{\mathbb{R}^n} \int_{\mathbb{R}^n} f(\mathbf{z}) 
d \mathbf{x} d \mathbf{y}$. Hence, 
\begin{eqnarray}
\mu_A (\mathbb{C}^n) 
&=& 2 \int_{\mathbb{C}^n} e^{-\frac{\mathbf{z}^T A \bar{\mathbf{z}}}{2}} 
d \mathbf{z} \nonumber\\
&=& 2 \int_{\mathbb{R}^n} \int_{\mathbb{R}^n} e^{-\frac{(\mathbf{x} + 
i \mathbf{y})^T (C + B) (\mathbf{x} - i \mathbf{y})}{2}} d \mathbf{x} d \mathbf{y} 
\nonumber\\
&=& 2 \int_{\mathbb{R}^n} \int_{\mathbb{R}^n} e^{-\frac{\mathbf{x}^T C \mathbf{x} + 
\mathbf{y}^T C \mathbf{y} - 2i \mathbf{x}^T B \mathbf{y}}{2}} d \mathbf{x} 
d \mathbf{y} \label{prptysksmc}\\
&=& 2 \int_{\mathbb{R}^n} e^{- \left( \frac{\mathbf{x}^T C \mathbf{x} + 
\mathbf{x}^T B C^{-1} B^T \mathbf{x}}{2} \right) } \left( \int_{\mathbb{R}^n} e^{- 
\left( \frac{(\mathbf{y} - iC^{-1} B^T \mathbf{x})^T C (\mathbf{y} - iC^{-1} B^T 
\mathbf{x})}{2} \right) } d \mathbf{y} \right) d \mathbf{x} \nonumber\\
&=& 2 \int_{\mathbb{R}^n} e^{- \left( \frac{\mathbf{x}^T C \mathbf{x} + 
\mathbf{x}^T B C^{-1} B^T \mathbf{x}}{2} \right) } \left( \sqrt{2 \pi} \right)^n 
|C|^{-\frac{1}{2}} d \mathbf{x} \label{prptymltnl}\\
&=& 2 \left( \sqrt{2 \pi} \right)^n |C|^{-\frac{1}{2}} \int_{\mathbb{R}^n} e^{- 
\frac{\mathbf{x}^T (C + BC^{-1} B^T) \mathbf{x}}{2}} d \mathbf{x} \nonumber\\
&=& 2 \left( \sqrt{2 \pi} \right)^{2n} |C|^{-\frac{1}{2}} \left| C + BC^{-1} B^T 
\right|^{-\frac{1}{2}} \nonumber\\
&=& 2 \left( \sqrt{2 \pi} \right)^{2n} |C|^{-\frac{1}{2}} |C|^{-\frac{1}{2}} 
\left| I - C^{-\frac{1}{2}} B C^{-\frac{1}{2}} C^{-\frac{1}{2}} B C^{-\frac{1}{2}} 
\right|^{-\frac{1}{2}} \nonumber\\
&=& 2 \left( \sqrt{2 \pi} \right)^{2n} |C|^{-\frac{1}{2}} |C|^{-\frac{1}{2}} \left| 
I - C^{-\frac{1}{2}} B C^{-\frac{1}{2}} \right|^{-\frac{1}{2}} \left| I + 
C^{-\frac{1}{2}} B C^{-\frac{1}{2}} \right|^{-\frac{1}{2}} \nonumber 
\end{eqnarray}

\noindent
Note that (\ref{prptysksmc}) follows from the fact that $C$ is symmetric and $B$ is 
skew-symmetric, and (\ref{prptymltnl}) follows by using properties of the 
multivariate normal distribution. Continuing further, 
\begin{eqnarray*}
\mu_A (\mathbb{C}^n) 
&=& 2 \left( \sqrt{2 \pi} \right)^{2n} |C - B|^{-\frac{1}{2}} 
|C + B|^{-\frac{1}{2}}\\
&=& 2 \left( \sqrt{2 \pi} \right)^{2n} \left| A^T \right|^{-\frac{1}{2}} 
|A|^{-\frac{1}{2}}\\
&=& 2 \left( \sqrt{2 \pi} \right)^{2n} |A|^{-1} 
\end{eqnarray*}

\noindent
If we define the measure $\mu_A^*$ on $\mathbb{C}^n$ by 
$$
\mu_A^* (.) = \frac{\mu_A (.)}{2 \left( \sqrt{2 \pi} \right)^{2n} |A|^{-1}}, 
$$

\noindent
then $\mu_A^* (\mathbb{C}^n) = 1$.

\smallskip

\noindent
Let $D$ denote a diagonal matrix with non-negative entries $d_1, d_2, ... ,d_n$. 
Then, 
\begin{eqnarray*}
\int_{\mathbb{C}^n} e^{- \frac{1}{2} \sum_{j=1}^n d_j z_j \bar{z}_j} d \mu_A^* 
(\mathbf{z}) 
&=& \int_{\mathbb{C}^n} \frac{d \mu_{D+A} (\mathbf{z})}{2 \left( \sqrt{2 \pi} 
\right)^{2n} |A|^{-1}}\\
&=& \frac{\mu_{D+A} (\mathbb{C}^n)}{2 \left( \sqrt{2 \pi} \right)^{2n} |A|^{-1}}\\
&=& \frac{2 \left( \sqrt{2 \pi} \right)^{2n} |D + A|^{-1}}{2 \left( \sqrt{2 \pi} 
\right)^{2n} |A|^{-1}} 
\end{eqnarray*}

\noindent
Hence, 
\begin{eqnarray*}
\int_{\mathbb{C}^n} e^{-\frac{1}{2} \sum_{j=1}^n d_j |z_j|^2} d \mu_A^* 
(\mathbf{z}) = \frac{|A|}{|D + A|}. 
\end{eqnarray*}

\noindent
Let $\mu_A^{abs}$ denote the measure induced on $\left( \frac{|z_1|^2}{2}, 
\frac{|z_2|^2}{2}, ... ,\frac{|z_n|^2}{2} \right)$ by $\mu_A^*$. Since 
$\mu_A^{abs}$ is a signed measure, we can break it up as $\mu_A^{abs} = 
\mu_A^{abs,+} - \mu_A^{abs,-}$, where $\mu_A^{abs,+}$ and $\mu_A^{abs,-}$ are 
non-negative measures. One wonders for which $A, \; \mu_A^{abs,-} \equiv 0$, so 
that $\mu_A^{abs}$ is a probability measure. If $A$ is the upper $n \times n$ 
submatrix of the negative of the generator of a strictly skip-free Markov process, 
then it follows by (\ref{jtlpctrnfm}) and the uniqueness of Laplace transform for a 
bounded measure, that $\mu_A^{abs}$ is a probability measure. 

\section{The Markov property}
We now undertake the task of ascertaining the conditions under which $\left\{ 
l_{T_n}^i \right\}_{0 \leq i \leq n-1}$ satisifes the Markov property for $n \geq 
3$. 
\newtheorem{lemma}{Lemma}
\begin{lemma}
If $Q_n$ (with $Q$ as in (\ref{asmpnskpfr})) is tridiagonal, then $\left\{ 
l_{T_n}^i \right\}_{0 \leq i \leq n-1}$ satisfies the Markov property. 
\end{lemma}

\noindent
{\it Proof} Since $Q_n$ is tridiagonal, we can obtain a symmetric matrix $Q_n^*$ 
from $Q_n$ by diagonal conjugation. It follows that 
$$
\left\{ l_{T_n}^i \right\}_{0 \leq i \leq n-1} \stackrel{d}{=} \frac{1}{2} 
(\mathbf{\eta}^2 + \tilde{\mathbf{\eta}}^2), 
$$

\noindent
where, $\mathbf{\eta}, \tilde{\mathbf{\eta}} \sim MVN_n \left( \mathbf{0}, 
(-Q_n^*)^{-1} \right)$ and $\mathbf{\eta}, \tilde{\mathbf{\eta}}$ are independent. 
By the standard properties of Gaussian Markov random fields developed in 
\cite{ruhdgmrfbk}, we get that 
\begin{equation} \label{etamrkvpry}
\eta_k \mid (\eta_{k-1}, \eta_{k-2}, ... ,\eta_1) \stackrel{d}{=} \eta_k \mid 
\eta_{k-1} \; \forall 2 \leq k \leq n. 
\end{equation}

\noindent
Let $t > 0$ be fixed. Then, 
\begin{eqnarray*}
& & {\bf{E}} \left[ e^{-t(\eta_k^2 + \tilde{\eta}_k^2)} \mid \eta_{k-1}^2 + 
\tilde{\eta}_{k-1}^2, ... , \eta_1^2 + \tilde{\eta}_1^2 \right]\\
&=& {\bf{E}} \left[ {\bf{E}} \left[ e^{-t(\eta_k^2 + \tilde{\eta}_k^2)} \mid 
\eta_{k-1}, \eta_{k-2}, ... , \eta_1, \tilde{\eta}_k, ... , \tilde{\eta}_1 \right] 
\mid \eta_{k-1}^2 + \tilde{\eta}_{k-1}^2, ... , \eta_1^2 + \tilde{\eta}_1^2 
\right]\\
&=& {\bf{E}} \left[ e^{-t \tilde{\eta}_k^2} {\bf{E}} \left[ e^{-t \eta_k^2} \mid 
\eta_{k-1}, \eta_{k-2}, ... , \eta_1 \right] \mid \eta_{k-1}^2 + 
\tilde{\eta}_{k-1}^2, ... , \eta_1^2 + \tilde{\eta}_1^2 \right] 
\end{eqnarray*}

\noindent
The previous equality follows by the independence of $\mathbf{\eta}$ and 
$\tilde{\mathbf{\eta}}$. By (\ref{etamrkvpry}) we get, 
\begin{eqnarray*}
& & {\bf{E}} \left[ e^{-t(\eta_k^2 + \tilde{\eta}_k^2)} \mid \eta_{k-1}^2 + 
\tilde{\eta}_{k-1}^2, ... , \eta_1^2 + \tilde{\eta}_1^2 \right]\\
&=& {\bf{E}} \left[ {\bf{E}} \left[ e^{-t \eta_k^2} \mid \eta_{k-1} \right] 
e^{-t \tilde{\eta}_k^2} \mid \eta_{k-1}^2 + \tilde{\eta}_{k-1}^2, ... , 
\eta_1^2 + \tilde{\eta}_1^2 \right]\\
&=& {\bf{E}} \left[ {\bf{E}} \left[ e^{-t \eta_k^2} \mid \eta_{k-1} \right] 
{\bf{E}} \left[ e^{-t \tilde{\eta}_k^2} \mid \tilde{\eta}_{k-1}, 
\tilde{\eta}_{k-2}, ... , \tilde{\eta}_1, \eta_{k-1}, ... , \eta_1 \right] \mid 
\eta_{k-1}^2 + \tilde{\eta}_{k-1}^2, ... , \eta_1^2 + \tilde{\eta}_1^2 \right]\\
&=& {\bf{E}} \left[ {\bf{E}} \left[ e^{-t \eta_k^2} \mid \eta_{k-1} \right] 
{\bf{E}} \left[ e^{-t \tilde{\eta}_k^2} \mid \tilde{\eta}_{k-1} \right] \mid 
\eta_{k-1}^2 + \tilde{\eta}_{k-1}^2, ... , \eta_1^2 + \tilde{\eta}_1^2 \right] 
\end{eqnarray*}

\noindent
The previous equality follows by the fact that $\mathbf{\eta}$ and 
$\tilde{\mathbf{\eta}}$ are i.i.d. Note that, 
$$
\eta_k \mid \eta_{k-1} \sim N(c_1 \eta_{k-1}, c_2) \mbox { and } 
\tilde{\eta}_k \mid \tilde{\eta}_{k-1} \sim N(c_1 \tilde{\eta}_{k-1}, c_2) 
$$

\noindent
for some constants $c_1$ and $c_2$. This gives, 
\begin{equation} \label{cndtliylpt}
{\bf{E}} \left[ e^{-t \eta_k^2} \mid \eta_{k-1} \right] = \frac{e^{-\frac{c_1^2 t 
\eta_{k-1}^2}{2c_2 t + 1}}}{\sqrt{2c_2 t + 1}} \mbox{ and } {\bf{E}} \left[ e^{-t 
\tilde{\eta}_k^2} \mid \tilde{\eta}_{k-1} \right] = \frac{e^{-\frac{c_1^2 t 
\tilde{\eta}_{k-1}^2}{2c_2 t + 1}}}{\sqrt{2c_2 t + 1}}. 
\end{equation}

\noindent
Combining everything, 
\begin{eqnarray*}
{\bf{E}} \left[ e^{-t(\eta_k^2 + \tilde{\eta}_k^2)} \mid \eta_{k-1}^2 + 
\tilde{\eta}_{k-1}^2, ... , \eta_1^2 + \tilde{\eta}_1^2 \right] 
&=& {\bf{E}} \left[ \frac{e^{-\frac{c_1^2 t (\eta_{k-1}^2 + \tilde{\eta}_{k-1}^2)}
{2c_2 t + 1}}}{2c_2 t + 1} \mid \eta_{k-1}^2 + \tilde{\eta}_{k-1}^2, ... , 
\eta_1^2 + \tilde{\eta}_1^2 \right]\\
&=& \frac{e^{-\frac{c_1^2 t (\eta_{k-1}^2 + \tilde{\eta}_{k-1}^2)}{2c_2 t + 1}}}
{2c_2 t + 1} 
\end{eqnarray*}

\noindent
Thus, the conditional Laplace transform given $\eta_{k-1}^2 + 
\tilde{\eta}_{k-1}^2, ... , \eta_1^2 + \tilde{\eta}_1^2$ is a function of 
$\eta_{k-1}^2 + \tilde{\eta}_{k-1}^2$. Since, 
\begin{eqnarray*}
\left\{ l_{T_n}^i \right\}_{0 \leq i \leq n-1} \stackrel{d}{=} \frac{1}{2} 
(\mathbf{\eta}^2 + \tilde{\mathbf{\eta}}^2), 
\end{eqnarray*}

\noindent
it follows that 
\begin{eqnarray*}
l_{T_n}^k \mid (l_{T_n}^{k-1}, l_{T_n}^{k-2}, ... , l_{T_n}^0) \stackrel{d}{=} 
l_{T_n}^k \mid l_{T_n}^{k-1} \; \; \forall 1 \leq k \leq n-1. 
\end{eqnarray*}

\noindent
Hence proved. 

\smallskip

\noindent
We now prove the converse. 
\begin{lemma}
If $Q_n$ (with $Q$ as in (\ref{asmpnskpfr})) is not tridiagonal, $\left\{ l_{T_n}^i 
\right\}_{0 \leq i \leq n-1}$ does not satisfy the Markov property. 
\end{lemma}

\noindent
{\it Proof} We first look at the scenario where $(X_1, X_2, X_3)$ are non-negative 
random variables with joint Laplace transorm given by 
$$
{\bf{E}} \left[ e^{-(d_1 X_1 + d_2 X_2 + d_3 X_3)} \right] = \frac{|-A|}{|D - A|}, 
$$

\noindent
where, 
$$
A = \left( \matrix{-a_{11} & a_{12} & 0 \cr a_{21} & -a_{22} & a_{23} \cr a_{31} & 
a_{32} & -a_{33}} \right), 
$$

\noindent
with $a_{12}, a_{23}, a_{31}, a_{22}, a_{33} > 0, \; a_{21}, a_{32} \geq 0$ and 
row sums less than zero or equal to zero. It follows that, 
$$
{\bf{E}} \left[ e^{-(d_2 X_2 + d_3 X_3)} \right] = \frac{|-A|}{a_{11} \left| 
\matrix{d_2 + \tilde{a}_{22} & -\tilde{a}_{23} \cr -\tilde{a}_{32} & d_3 + a_{33}} 
\right|}, 
$$

\noindent
where $\tilde{a}_{22} = a_{22} - \frac{a_{12} a_{21}}{a_{11}}, \; \tilde{a}_{23} = 
\tilde{a}_{32} = \sqrt{a_{23} a_{32} + \frac{a_{12} a_{23} a_{31}}{a_{11}}}$. It 
follows that, 
$$
(X_2, X_3) \stackrel{d}{=} \frac{1}{2} (\eta_1^2, \eta_2^2) + \frac{1}{2} 
(\tilde{\eta}_1^2, \tilde{\eta}_2^2), 
$$

\noindent
where $\mathbf{\eta}, \tilde{\mathbf{\eta}} \sim MVN_2 \left( \mathbf{0}, 
(-A^*)^{-1} \right)$ are independent with 
$$
A^* = \left( \matrix{-\tilde{a}_{22} & \tilde{a}_{23} \cr \tilde{a}_{32} & 
-a_{33}} \right). 
$$

\noindent
By (\ref{cndtliylpt}), 
$$
{\bf{E}} \left[ e^{-d_3 X_3} \mid X_2 \right] = \frac{a_{33}}{d_3 + a_{33}} 
e^{-\frac{(\tilde{a}_{23})^2 d_3 X_2}{a_{33} (d_3 + a_{33})}}. 
$$

\noindent
Note that, 
\begin{eqnarray*}
& & X_3 \mid X_2, X_1 \stackrel{d}{=} X_3 \mid X_2\\
&\Leftrightarrow& {\bf{E}} \left[ {\bf{E}} \left[ e^{-d_3 X_3} \mid X_2, X_1 
\right] e^{-(d_1 X_1 + d_2 X_2)} \right] = {\bf{E}} \left[ {\bf{E}} \left[ e^{-d_3 
X_3} \mid X_2 \right] e^{-(d_1 X_1 + d_2 X_2)} \right] \; \forall d_1, d_2, d_3 
\geq 0\\
&\Leftrightarrow& \frac{|-A|}{|D - A|} = {\bf{E}} \left[ \frac{a_{33}}{d_3 + 
a_{33}} e^{-(d_1 X_1 + d_2^* X_2)} \right], \mbox{ with } d_2^* = d_2 + 
\frac{(\tilde{a}_{23})^2 d_3}{a_{33} (d_3 + a_{33})} \; \forall d_1, d_2, d_3 
\geq 0\\
&\Leftrightarrow& \frac{|-A|}{|D - A|} = \frac{a_{33}}{d_3 + a_{33}} \frac{|-A|}
{(d_1 + a_{11})(d_2^* + a_{22})a_{33} - (d_1 + a_{11})a_{23}a_{32} - a_{12} a_{21} 
a_{33} - a_{12} a_{23} a_{31}}\\
& & \forall d_1, d_2, d_3 \geq 0 
\end{eqnarray*}

\noindent
After expanding the determinant and cancelling terms from both sides we obtain, 
\begin{eqnarray*}
& & X_3 \mid X_2, X_1 \stackrel{d}{=} X_3 \mid X_2\\
&\Leftrightarrow& \frac{(d_1 + a_{11})a_{12} a_{23} a_{31} d_3}{a_{33} a_{11}} - 
\frac{(d_3 + a_{33}) a_{12} a_{23} a_{31}}{a_{33}} = 0 \; \forall d_1, d_2, d_3 
\geq 0\\
&\Leftrightarrow& a_{12} a_{23} a_{31} (d_1 d_3 - a_{11}a_{33}) = 0 \; \forall d_1, 
d_2, d_3 \geq 0, 
\end{eqnarray*}

\noindent
which is not true, because $a_{12}, a_{23}, a_{31} > 0$ by our assumptions. Hence, 
$(X_1, X_2, X_3)$ is not Markov. 

\indent
We now return to $\left\{ l_{T_n}^i \right\}_{0 \leq i \leq n-1}$. The scenario 
considered above immediately gives the required result for $n = 3$. But it also 
helps us in the general case. Since $Q_n$ is not tridiagonal, we can find the 
smallest $i$ such that $q_{ij} > 0$ for some $j \leq i-2$. Let us call it as $i_0$ 
and evaluate the joint Laplace transform of $\left( l_{T_n}^{i_0 - 2}, 
l_{T_n}^{i_0 - 1}, l_{T_n}^{i_0} \right)$. By (\ref{jtlpctrnfm}), 
$$
{\bf{E}} \left[ e^{-\left( d_1 l_{T_n}^{i_0 - 2} + d_2 l_{T_n}^{i_0 - 1} + d_3 
l_{T_n}^{i_0} \right)} \right] = \frac{|-Q_n|}{\left| D^{i_0,3} - Q_n \right|}, 
$$

\noindent
where $D^{i_0,3}$ is a diagonal matrix with the $(i_0 - 2)^{th}, (i_0 - 1)^{th}$ 
and $i_0^{th}$ entries given by $d_1, d_2$ and $d_3$ respectively and other 
diagonal entries being $0$. We now reduce the determinant $\left| D^{i_0,3} - Q_n 
\right|$ to a managable form by performing appropriate row and column operations. 
Subtract a multiple of the first column from the second column to make the 
$(1,2)^{th}$ entry zero. Successively subtract a multiple of the $(j-1)^{th}$ 
column from the $j^{th}$ column to make the $(j-1,j)^{th}$ entry $0$, for $j \leq 
i_0 - 2$. The upper $(i_0 - 3) \times (i_0 - 3)$ submatrix is now upper triangular. 
Similarily, subtract a multiple of the $n^{th}$ row from the $(n-1)^{th}$ row to 
make the $(n-1,n)^{th}$ entry $0$. Successively subtract a multiple of the $j^{th}$ 
row from the $(j-1)^{th}$ row to make the $(j-1,j)^{th}$ entry $0$, for $j \geq i_0 
+ 1$. The lower $(n - i_0) \times (n - i_0)$ matrix is now upper triangular. The 
fact that $Q$ is a generator matrix as in (\ref{asmpnskpfr}) ensures that all the 
diagonal entries are positive. Hence, 
$$
\left| D^{i_0,3} - Q_n \right| = c \left| \matrix{d_1 - \tilde{q}_{i_0 - 2, 
i_0 - 2} & -q_{i_0 - 2, i_0 - 1} & 0 \cr -q_{i_0 - 1, i_0 - 2} & d_2 - q_{i_0 - 1, 
i_0 - 1} & -q_{i_0 - 1, i_0} \cr -\tilde{q}_{i_0, i_0 - 2} & -\tilde{q}_{i_0, i_0 - 
1} & d_3 - \tilde{q}_{i_0, i_0}} \right| \mbox{ with } c > 0. 
$$

\noindent
Again, $Q$ being a generator matrix as in (\ref{asmpnskpfr}) along with $q_{i_0 j} 
> 0$ for some $j \leq i_0 - 2$ ensures that 
$$
Q_n^{i,3} = \left( \matrix{\tilde{q}_{i_0 - 2, i_0 - 2} & q_{i_0 - 2, i_0 - 1} & 0 
\cr q_{i_0 - 1, i_0 - 2} & q_{i_0 - 1, i_0 - 1} & q_{i_0 - 1, i_0} \cr 
\tilde{q}_{i_0, i_0 - 2} & \tilde{q}_{i_0, i_0 - 1} & \tilde{q}_{i_0, i_0}} \right) 
$$

\noindent
is a generator matrix with $q_{i_0 - 2, i_0 - 1}, q_{i_0 - 1, i_0}, \tilde{q}_{i_0, 
i_0 - 2} > 0$ and row sums less than zero or equal to zero. Hence, 
$$
{\bf{E}} \left[ e^{-\left( d_1 l_{T_n}^{i_0 - 2} + d_2 l_{T_n}^{i_0 - 1} + d_3 
l_{T_n}^{i_0} \right)} \right] = \frac{\left| - Q_n^{i_0,3} \right|}{\left| 
D^3 - Q_n^{i_0,3} \right|}, 
$$

\noindent
where $D^3$ is a diagonal matrix of dimension $3$ with diagonal entries $d_1, d_2, 
d_3$. The Laplace transform of $ \left( l_{T_n}^{i_0 - 2}, l_{T_n}^{i_0 - 1}, 
l_{T_n}^{i_0} \right)$ satisfies the conditions of the scenario considered at the 
beginning of this proof. Hence, $\left( l_{T_n}^{i_0 - 2}, l_{T_n}^{i_0 - 1}, 
l_{T_n}^{i_0} \right)$ is not Markov, which is enough to conclude that $\left\{ 
l_{T_n}^i \right\}_{0 \leq i \leq n-1}$ is not Markov. This completes the proof.


\begin{thebibliography}{5}
\bibitem{dynkmrkvps} Dynkin, E.B. (1983). Markov processes as a tool in field 
theory, {\it{Journal of Functional Analysis}} {\bf{50}}, 167-187. 
\bibitem{eikamarosh} Eisenbaum, N., Kaspi, H., Marcus, M., Rosen, J. and Shi, Z. 
(2000). A Ray Knight theorem for symmetric Markov processes, {\it{Annals of 
Probability}} {\bf{28}}, 1781-1796. 
\bibitem{kentmedocc} Kent, J.T. (1983). The appearance of a multivariate 
exponential distribution in sojourn times for birth-death and diffusion processes, 
{\it{Probability, Statistics and Analysis, Edited by J.F.C. Kingman and 
G.E.H. Reuter, London Mathematical Society Lecture Note Series}} {\bf{79}}, 
Cambridge University Press. 
\bibitem{lejandknim} Le Jan, Y. (2006). Dynkin's isomorphism without symmetry, to 
appear in {\it{Proceedings of the conference ``Stochastic Analysis in Mathematical 
Physics"}}, Lisbon 2006. 
\bibitem{ruhdgmrfbk} Rue, H. and Held, L. (2005). {\it{Gaussian Markov Random 
Fields: Theory and Applications (Monographs on Statistics and Applied 
Probability)}}, Chapman and Hall, Boca Raton, Florida. 
\end{thebibliography}
\end{document}